\documentclass[11pt,a4paper]{article} 
\usepackage[cp1251]{inputenc}         
\usepackage{amssymb,amsfonts,amsmath} 
\usepackage{graphicx}                 
\usepackage{enumerate}                %

\newtheorem{pr}{Proposition}[section]
\newtheorem{cor}{Corollary}[section]
\newtheorem{lm}{Lemma}[section]
\newtheorem{rem}{Remark}[section]

\newtheorem{exam}{Example}[section]

\begin{document}

\title{On quotients of affine superschemes over finite supergroups}
\author{A.N.Zubkov}
\date{}

\maketitle

\begin{abstract}
In this article we consider sheaf quotients of affine superschemes
by finite supergroups that act on them freely. More precisely, if
a finite supergroup $G$ acts on an affine superscheme $X$ freely,
then the quotient $K$-sheaf $\tilde{X/G}$ is again an affine
superscheme $Y$, where $K[Y]\simeq K[X]^G$. Besides, $K[X]$ is a
finitely presented projective $K[X]^G$-module.
\end{abstract}

\section*{Introduction}

In the present article we prove that if a finite supergroup $G$
acts on an affine superscheme $X$ freely, then the sheaf quotient
$\tilde{X/G}$ is again affine and isomorphic to $SSp \ R$, where
$R=K[X]^G$. Moreover, we also prove that $K[X]$ is finitely
presented projective $R$-module. This theorem generalizes the
classical, purely even case (cf. \cite{dg, jan}). On the whole, we
follow the ideas from \cite{dg} but there is a principal
difference between purely even and super cases.  In the classical
case $K[X]$ is always integral over $R$. In the supercase it is
not still true (see Example 3.1 below)! It happens as soon as the
$G$-action is not free. In an equivalent formulation, for some
finite supergroups 14-th Hilbert problem has the negative
solution. To overcome this obstacle we exploit the freeness of our
action and reduce the general case to the case, when $G$ has not
any proper normal supersubgroups.

\section{Superalgebras and supermodules}

In what follows all superalgebras are commutative. The category of
commutative superalgebras with even morphisms is denoted by
$SAlg_K$. If $A\in SAlg_K$, then the category of left (right)
$A$-supermodules with even morphisms is denoted by $A-smod$
(respectively, $smod-A$). Remind that $A-smod\simeq smod-A$
\cite{z}. More precisely, any $M\in A-smod$ has the structure of a
right $A$-supermodule via $ma=(-1)^{|a||m|}am, a\in A, m\in M$.

Remind that any left or right maximal ideal ${\cal M}$ of a
superalgebra $A$ is a two-sided superideal \cite{z}, Lemma 1.1.
Moreover, ${\cal M}={\cal M}_0\bigoplus A_1$, where ${\cal M}_0$
is a maximal ideal of $A_0$. Let $M$ be a free $A$-supermodule of
(finite) superrank $(m, n)$. Take elements $m_1 ,\ldots ,
m_{m+n}\in M$ such that $|m_i|=0$ iff $1\leq i\leq m$, otherwise
$|m_i|=1$.
\begin{lm} The elements $m_1 ,\ldots , m_{m+n}$ form a free
basis of $M$ iff their canonical images form a free basis of
$A/rad A$-supermodule $M/(rad A)M$.
\end{lm}
Proof. By Nakayama's lemma the supersubmodule $N=\sum_{1\leq i\leq
m+n} Am_i$ coincides with $M$ (cf. \cite{kash}, Theorem 9.2.1(d)).
Lemma 5.5 from \cite{z} concludes the proof.

We say that an $A$-supermodule $M$ is finitely generated, if $M$
is an epimorphic image of a free $A$-supermodule of finite
superrank. Besides, if the kernel of the above epimorphism is also
finitely generated, then $M$ is called {\it finitely presented}.
It is obvious that $M$ is finitely generated as a supermodule iff
it is finitely generated as a module.
\begin{lm}
A supermodule $M$ is finitely presented iff it is finitely
presented as an $A$-module.
\end{lm}
Proof. Let $\phi : A^n=\bigoplus_{1\leq i\leq n} A e_i\to M$ be an
epimorphism of $A$-modules such that $\ker\phi$ is a finitely
generated $A$-submodule of $A^n$. Consider a free $A$-supermodule
$A^{n|n}$ with a basis $e_{i, \epsilon}, |e_{i,
\epsilon}|=\epsilon, 1\leq i\leq n, \epsilon=0, 1$. Denote
$\phi(e_i)$ by $m_i$. Define the supermodule epimorphism $\psi :
A^{n|n}\to M$ by $\psi(e_{i, \epsilon})=m_{i, \epsilon}$, where
$|m_{i, \epsilon}|=\epsilon$ and $m_{i, 0}+ m_{i, 1}=m$. Since the
elements $m_i$ generate $M$, we have $m_{i, 0}=\sum_{1\leq j\leq
n}a_{ij}m_j, 1\leq i\leq n, a_{ij}\in A$. A diagram
$$\begin{array}{lcr}
A^n & \stackrel{\phi}{\to} & M \\
p \nwarrow & &\nearrow\psi   \\
 & A^{n|n} &
\end{array}, $$
where $p(e_{i, 0})=\sum_{1\leq j\leq n}a_{ij}e_j, p(e_{i,
1})=e_i-p(e_{i, 0})$, is obviously commutative. Moreover, $p$ is
an epimorphism and $\ker p$ contains a submodule $T$, generated by
the elements $$e_{i, 0}-\sum_{1\leq j\leq n}a_{ij}(e_{j, 0}+e_{j,
1}), 1\leq i\leq n .$$ As $A^{n|n}/T$ is generated by the residue
classes of $n$ elements $(e_{i, 0}+e_{i, 1})$, it follows that $p$
induces an isomorphism $A^{n|n}/T\simeq A^n$. In particular, the
supersubmodule $\ker\psi=p^{-1}(\ker\phi)$ is finitely generated.
\begin{rem}
If a superalgebra $A$ is finitely presented as a module over its
supersubalgebra $B$, then $A$ is finitely presented as a
$B$-superalgebra.
\end{rem}
A superalgebra $A$ is called semi-local iff $A$ contains only
finitely many maximal ideals. By the above, $A$ is semi-local iff
$A_0$ is semi-local. Let ${\cal N}_1,\ldots , {\cal N}_t$ are all
maximal ideals of $A$. It can be easily checked that {\it Chinese
reminder Theorem} holds for two-sided ideals of any (not necessary
commutative) algebra or ring (see for example \cite{bur}, II, \S
1, Proposition 5). Thus
$$A/ rad A=A/\bigcap_{1\leq i\leq t}{\cal N}_i\simeq\prod_{1\leq i\leq t}
A/ {\cal N}_i$$ is a direct product of fields. Conversely, if $A/
rad A$ is a direct product of finitely many fields, then $A$ is
semi-local. Besides, if $L$ is an $A$-module, then
$$L/(rad A)L\simeq\prod_{1\leq i\leq t}L/{\cal N}_i L .$$

Let $A$ be a semi-local superalgebra and $B$ be its local
supersubalgebra whose maximal ideal ${\cal M}$ is contained in
$rad A$. Let $M$ be a free $A$-supermodule of finite superrank.
\begin{lm}
If $N$ is a $B$-supersubmodule of $M$ such that $AN=M$ and
$B/{\cal M}$ is an infinte field, then $N$ contains a free basis
of $M$.
\end{lm}
Proof. Using Lemma 1.1 one can replace $A, M, N$ by $A/ rad A ,
M/(rad A)M$ and $(N+ (rad A)M)/(rad A)M$ respectively. The final
arguing can be copied from \cite{dg}, III, \S 2, Lemma 4.7.

Let $C\in SAlg_K$ and $S$ is a multiplicative subset of $C_0$. One
can define its left (right) $S^{-1}C$-supermodule of fractions
$S^{-1}M=S^{-1}C\otimes_C M$ (respectively, $M S^{-1}=M\otimes_C
S^{-1}C$). It is clear that $S^{-1} M$, considered as a right
$S^{-1}C$-supermodule, is isomorphic to $M S^{-1}$. The
isomorphism is given by $x\otimes m\to (-1)^{|x||m|}m\otimes x,
x\in S^{-1}C, m\in M$.
\begin{lm}
If $M, N\in C-smod$, then $S^{-1}C$-supermodules
$$S^{-1}M\otimes_C N , M\otimes _C S^{-1}N,
S^{-1}M\otimes_{S^{-1}C}S^{-1}N \ \mbox{and} \ S^{-1}(M\otimes_C
N)$$ are canonically isomorphic each to other.
\end{lm}
Proof. Routine checking (see also \cite{bur}, II, \S 2,
Proposition 18).

Let $A$ be a superalgebra and $B$ be its supersubalgebra. We say
that $A$ is an {\it integral extension} of $B$ (or $A$ is {\it
integral} over $B$) iff $A_0$ is an integral extension of $B_0$.
The following lemma is an obvious consequence of Lemma 1.2,
\cite{z}.
\begin{lm}
If $B\subseteq A$ is integral and $\cal P$ is a prime ideal of
$B$, then there is a prime ideal $\cal Q$ of $A$ such that ${\cal
Q}\bigcap B={\cal P}$ (in particular, ${\cal P}A\neq A$).
Moreover, $\cal P$ is maximal iff $\cal Q$ is maximal.
\end{lm}
We say that $\cal Q$ {\it lies over} $\cal P$. Notice that a
maximal ideal of the superalgebra of fractions $A_{\cal
P}=(B_0\setminus {\cal P}_0)^{-1}A$ has a form ${\cal Q}_{\cal
P}$, where $\cal Q$ lies over $\cal P$. It infers that $A_{\cal
P}$ is semi-local iff there are finitely many prime ideals of $A$
those lie over $\cal P$. The last property is guarantied for any
$\cal P$, whenever $A_0$ is a finitely generated $B_0$-module (cf.
\cite{bur}, V, \S 2, Proposition 3).

The proof of the following lemma can be copied from Proposition 8
and Proposition 9, \cite{bur}, II, \S 3.
\begin{lm}
Let $\cal M$ be a maximal ideal of a superalgebra $A$ and $M$ be
an $A$-supermodule. The canonical morphism $M/{\cal M}^t M\to
M_{\cal M}/{\cal M}^t M_{\cal M}$ is a supermodule isomorphism.
\end{lm}
\begin{lm}
Let $u: M\to N$ be a morphism of $A$-supermodules. If $N$ is
finitely generated, then $u$ is surjective iff for any maximal
ideal $\cal M$ of $A$ the induced morphism $M/{\cal M}M\to N/{\cal
M}N$ is surjective.
\end{lm}
Proof. Combine Lemma 1.6 with Lemma 1.5 from \cite{z} and argue as
in Proposition 11, \cite{bur}, II, \S 3.
\begin{rem}
The statements of the above lemmas are still true, even if $M$ and
$N$ are $A$-modules.
\end{rem}
\begin{pr}
Let $\phi : B\to A$ be a superalgebra morphism such that $A$ is
finitely generated $B$-module. If the induced morphism $SSp \ A\to
SSp \ B$ is an inclusion of $K$-functors, then $\phi$ is an
epimorphism.
\end{pr}
Proof. We use the following nice trick from \cite{dg}, I, \S 5,
1.5 . The diagonal morphism $SSp \ A\to SSp \ A\times_{SSp \ B}
SSp \ A$ is an isomorphism. By Yoneda's lemma the canonical
superalgebra morphism $A\otimes_B A\to A$ by $a_1\otimes a_2\to
a_1 a_2, a_1, a_2\in A,$ is an isomorphism. It implies that for
any maximal ideal $\cal N$ of $B$ we have the isomorphism $A/{\cal
N}A\otimes_{B/{\cal N}} A/{\cal N}A\to A/{\cal N}A$. Comparing
dimensions (over the field $B/{\cal N}$) we see that $B/{\cal
N}\to A/{\cal N}$ is surjective. Lemma 1.7 concludes the proof.

\section{Unipotent supergroups}

In what follows a supersubgroup of an affine or algebraic
supergroup is closed. We use notations and definitions from
\cite{z}.

Let $G$ be an algebraic supergroup. It is called {\it unipotent}
iff any simple $G$-supermodule is one dimensional and trivial. It
is easy to see that $G$ is unipotent iff for any non-zero
$G$-supermodule $V$ its invariant subspace $V^G$ is not zero also.
By Proposition 6.2 from \cite{z} there is a finite-dimensional
$G$-supermodule $V$ such that $G$ is isomorphic to an
supersubgroup of $GL(V)$. Since $G$ is unipotent, there is a flag
of $G$-supersubmodules
$$0\subseteq V_1\subseteq\ldots\subseteq V_r=V$$ such that for all
$i\geq 1$ $V_i/V_{i-1}$ is a trivial $G$-supermodule. Denote this
flag by $\bf V$. Consider the subfunctor $U({\bf V})\subseteq
GL(V)$ defined by $$U({\bf V})(A)=\{g\in Stab_{\bf V}(A)|
g|_{V_i\otimes A} \ \mbox {acts identically modulo} \
V_{i-1}\otimes A, i\geq 1\}, A\in SAlg_K.$$ It is clear that
$U({\bf V})$ is a supersubgroup of $GL(V)$. In fact, there is a
basis $v_1,\ldots , v_r$ of the superspace $V$ such that $|v_i|=0$
if $1\leq i\leq m$, otherwise $|v_i|=1$, and a unique substitution
$\sigma\in S_r$ with $\sigma(1)<\ldots \sigma(m),
\sigma(m+1)<\ldots\ldots <\sigma(r)$. Besides, $v_i$ generates
$V_{\sigma(i)}/V_{\sigma(i)-1}$. The supersubgroup $U({\bf V})$ is
defined by $x_{ii}=1$ and $x_{ji}=0, 1\leq i, j\leq r,
\sigma(j)>\sigma(i)$. We also use for $U({\bf V})$ the other
notations, say $U_{\sigma}$ or $U_{\sigma}(m|n)$, where $n=\dim
V_1=r-m$. By the above, $G\leq U_{\sigma}$.
\begin{rem}
A supergroup $U_{\sigma}$ is contained in $SL(V)$, that is
$Ber(U_{\sigma})=1$. In fact, the berezinian $Ber$  induces a
supergroup epimorphism $GL(V)\to GL(1|0)=G_m$. In particular,
$$K[U_{\sigma}]=K[x_{ij}|\sigma(i) <\sigma(j)]\simeq K[x_{ij}|
1\leq i < j\leq r]\simeq K[{\bf
A}^{\frac{m(m-1)}{2}+\frac{n(n-1)}{2}|mn}].$$
\end{rem}
\begin{lm}
Let $G$ be an algebraic supergroup. Assume that the superalgebra
$K[G]$ has a $G$-supermodule (or equivalently, a right
$K[G]$-supercomodule) filtration
$$0\subseteq W_1\subseteq W_2\subseteq\ldots $$ such that
$\bigcup_{i\geq 1}W_i=K[G]$ and each factor $W_i/W_{i-1}$ is a
trivial $G$-supermodule. Then $G$ is unipotent.
\end{lm}
Proof. Let $W$ be a simple $G$-supermodule and $f\in M^*, f\neq
0$. We have a supermodule morphism $g=(f\otimes id_{K[G]})\tau_W :
W\to K[G]$ of the same parity as $f$. Since the preimages
$g^{-1}(W_i)$ form a $G$-supermodule filtration of $W$ and $W$ is
simple, we see that $W$ is isomorphic to a factor $W_i/W_{i-1}$.
\begin{pr}
If $G$ is isomorphic to a supersubgroup of $U_{\sigma}(m|n)$, then
$G$ is unipotent.
\end{pr}
Proof. One has to build an $U_{\sigma}$-supermodule filtration of
$K[U_{\sigma}]$ as in Lemma 2.1. The superalgebra $R=K[U_{\sigma}]
$ has a natural $\bf N$-grading $R=\bigoplus_{k\geq 0} R_k$, where
each $R_k$ is a superspace. Ascribe to any monomial $m=x_{i_1
j_1}\ldots x_{i_k j_k}\in R_k\setminus 0$ the weight
$v(m)=\sum_{1\leq t\leq k}(\sigma(j_t)-\sigma(i_t))$. It is easy
to see that
$$\tau_{R}(m)-m\otimes 1\in\sum_{m'\in R_k, v(m')< v(m)} m'\otimes
R +\sum_{0\leq s\leq k-1} R_s\otimes R .$$ In particular, we have
a $U_{\sigma}$-supermodule filtration
$$0\subseteq K=R_{1, 0}\subseteq R_{1, 2}\subseteq\ldots\subseteq R_{k,
t}\subseteq\ldots ,$$ where
$$R_{k, t}=\bigoplus_{0\leq s\leq k-1}R_s\bigoplus (\sum_{m\in
R_k,\ v(m)\leq t}Km), \ k\geq 1, \ 0\leq t\leq (r-1)k .$$ Besides,
all sequential factors of this filtration are sums of trivial
$U_{\sigma}$-supermodules. Proposition is proved.
\begin{cor}
If $G$ is unipotent, then any its supersubgroup and
superfactorgroup is also unipotent.
\end{cor}

The supergroups $U_{1}(2|0)$ and $U_{1}(1|1)$ are usually denoted
by $G_a$ and $G_a^-$ respectively. Besides, $G_a$ is called {\it
even} and $G_a^-$ is called {\it odd} one-dimensional additive
supergroup.
\begin{lm}
Let $G$ be an algebraic supergroup and $N$ be its normal
supersubgroup such that $N$ and $\tilde{G/N}$ are unipotent. Then
$G$ is unipotent.
\end{lm}
Proof. Let $V$ be a simple $G$-supermodule. We know that $V^N\neq
0$ and $V^N$ is the largest supersubspace of $V$ whose coefficient
(super)space belongs to $K[G]^N=K[\tilde{G/N}]$. Thus $V=V^N$ is a
simple $G/N$-supermodule. In particular, $V$ is one-dimensional
and trivial.
\begin{lm}
Let $\pi : G\to H$ be an algebraic group epimorphism with the
kernel $N$. If $L$ is a supesubgroup of $G$, then
$\tilde{L/L\bigcap N}$ is canonically isomorphic to $Im\pi |_L$.
\end{lm}
Proof. Notice that $\ker\pi |_L=L\bigcap N$ and use Theorem 6.1
from \cite{z}.
\begin{pr}
If $G$ is unipotent, then $G$ has a series of normal
supersubgroups $1\leq N_1\unlhd\ldots\unlhd N_t=G$ such that any
factor $\tilde{N_i/N_{i-1}}$ belongs to $Z(\tilde{G/N_{i-1}})$ and
isomorphic either to a supersubgroup of $G_a$ or to $G_a^-$.
\end{pr}
Proof. By Lemma 2.3 all we have to prove is that such series
exists in $U_{\sigma}$. For any $k\geq 1$ define a superideal
$I_k$ of $K[U_{\sigma}]$, generated by the elements $x_{ij}$ with
$\sigma(j)-\sigma(i)\leq k$. It can easily be checked that
$U_{\sigma , k}=V(I_k)\unlhd U_{\sigma}$. Indeed, the superalgebra
$B_k=K[x_{ij}|\sigma(j)-\sigma(i)\leq k]$ is a Hopf
supersubalgebra of $K[U_{\sigma}]$ and $U_{\sigma , k}$ coincides
with the kernel of the epimorphism $U_{\sigma}\to SSp \ B_k$. In
the same way,
$$\tilde{U_{\sigma , k}/U_{\sigma , k+1}}\simeq SSp \
B_{k+1}/B_{k+1}B_k^+\simeq SSp \
K[x_{ij}|\sigma(j)-\sigma(i)=k+1]\simeq (G_a)^s\times (G_a^-)^l,$$
where $s$ (respectively, $l$) is the number of even (respectively,
odd) elements among $\{x_{ij}|\sigma(j)-\sigma(i)=k+1 \}$. It
remains to check that $\tilde{U_{\sigma , k}/U_{\sigma , k+1}}\leq
Z(\tilde{U_{\sigma}/U_{\sigma , k+1}})$. It is equivalent to the
statement that the superalgebra morphism
$$K[\tilde{U_{\sigma , k}/U_{\sigma , k+1}}]\otimes
K[\tilde{U_{\sigma}/U_{\sigma , k+1}}]\to K[\tilde{U_{\sigma ,
k}/U_{\sigma , k+1}}],$$ induced by $\nu_l$, coincides with the
morphism $f\to f\otimes 1, f\in K[\tilde{U_{\sigma , k}/U_{\sigma
, k+1}}]$. The last one is dual to the projection
$$\tilde{U_{\sigma , k}/U_{\sigma , k+1}}\times\tilde{U_{\sigma}/U_{\sigma ,
k+1}}\to \tilde{U_{\sigma , k}/U_{\sigma , k+1}}.$$ Since
$$\nu_l(x_{ij})=\sum_{t_1, t_2,
\sigma(i)\leq\sigma(t_1)<\sigma(t_2)\leq\sigma(j)} (-1)^{|x_{t_1 ,
t_2}||x_{i, t_1}|}x_{t_1, t_2}\otimes x_{i,
t_1}s_{U_{\sigma}}(x_{t_2, j}),
$$
we obtain that $\nu_l(x_{ij})=x_{ij}\otimes 1$ modulo
$B_{k+1}B_k^+$. Proposition is proved.

\section{Proof of the main theorem}

Let $G$ be an algebraic supergroup. Assume that $G$ acts on an
affine superscheme $X$. Denote the corresponding morphism of
(affine) superschemes $X\times G\to X$ by $\phi_0$. For the
reader's convenience we remind some basic notations and facts from
\cite{dg}, III, \S 2-4. The squares
$$\begin{array}{ccc}
X\times G\times G & \stackrel{\phi_0'}{\to} & X\times G \\
\phi'_2 \downarrow & & \downarrow\phi_1 \\
X\times G &\stackrel{\phi_0}{\to} & X
\end{array},
\begin{array}{ccc}
X\times G\times G & \stackrel{\phi_1'}{\to} & X\times G \\
\phi'_0 \downarrow & & \downarrow\phi_0 \\
X\times G &\stackrel{\phi_0}{\to} & X
\end{array},
\begin{array}{ccc}
X\times G\times G & \stackrel{\phi_1'}{\to} & X\times G \\
\phi'_2 \downarrow & & \downarrow\phi_1 \\
X\times G &\stackrel{\phi_1}{\to} & X
\end{array}
$$
are cartesian, where $\phi_1=pr_X, \phi'_2 =pr_{X\times G}$ and
$\phi'_0(x, g, h)=(xg, h), \phi'_1(x, g, h)=(x, gh), x\in X(A), g,
h\in G(A), A\in SAlg_K$. The morphism of superalgebras $K[X]\to
K[X]\otimes K[G]$, dual to $\phi_0$ (respectively, dual to
$\phi_1$), is denoted by $\tau_X$ (respectively, by $i_X$). The
supersubalgebra of (co)invariants $K[X]^G =\ker(\tau_X -i_X)$ is
denoted by $R$. Since $\phi_0$ has a left inverse $\sigma(x)=(x,
1), x\in X(A), A\in SAlg_k$, the couple $(X, \phi_0)$ is a
cokernel of the pair morphisms $(\phi'_0 , \phi'_1)$ (in the
category of $K$-functors!). Dualizing we obtain a commutative
diagram
$$
\begin{array}{ccccc}
K[X]\otimes K[G]^{\otimes 2} &
\begin{array}{c}{\stackrel{\delta'_0}{\leftarrow}} \\
{\stackrel{\delta'_1}{\leftarrow}}\end{array} & K[X]\otimes K[G]
&\stackrel{\tau_X}{\leftarrow} & K[X] \\
\delta'_2 \uparrow & & i_X \uparrow & & i \uparrow \\
K[X]\otimes K[G] & \begin{array}{c}{\stackrel{\tau_X}{\leftarrow}} \\
{\stackrel{i_X}{\leftarrow}}\end{array} & K[X]
&\stackrel{i}{\leftarrow} & R
\end{array},
$$
where $\delta'_0=\tau_X\otimes id_{K[G]}, \delta'_1=
id_{K[X]}\otimes\delta_G, \delta'_2=id_{K[X]\otimes K[G]}\otimes
1$. Its horizontal lines are exact and the left square is composed
from cocartesian squares those are dual to the above first and
third cartesian ones. We call this diagram {\it basic}.

From now on we assume that all supergroups are finite unless
otherwise stated. Without loss of generality one can assume that
$K$ is algebraically closed. The $K$-functor morphism $(\phi_1,
\phi_0) : X\times G\to X\times_{SSp \ R} X$ is dual to the
morphism of superalgebras
$$\psi : K[X]\otimes_R K[X]\to K[X]\otimes K[G]$$ defined as
$$f\otimes h\to \sum fh_1\otimes h_2, \tau_X(h)=\sum h_1\otimes
h_2, f, h, h_1\in K[X], h_2\in K[G]$$ (see \cite{z, jan}). We
denote $\dim K[G]$ by $|G|$ and call it the {\it order} of $G$.
The maximal ideal $\ker\epsilon_G$ is denoted by ${\cal M}$.
\begin{lm}
If $G$ acts on $X$ freely, then $\psi$ is surjective.
\end{lm}
Proof. Notice that $(\phi_1, \phi_0)$ is an injective $K$-functor
morphism and $Im\psi$ contains $K[X]\otimes 1$. It remains to
refer to Proposition 1.1.
\begin{lm}
Let $B$ be a supersubalgebra of a superalgebra $A$. Then : \\
1) If $A$ is a finitely generated $B$-module, then $A$ is integral over $B$; \\
2) If $A$ is a finitely generated superalgebra and integral over
$B$, then $A$ is a finitely generated $B$-module and $B$ is a
finitely generated superalgebra.
\end{lm}
Proof. To prove the first statement we fix a finite set of
generators of $B_0$-module $A_0/B_1 A_1$. Using Cayley-Hamilton's
theorem we see that for any $a\in A_0$ there is a unitary
polynomial $f(t)\in B_0[t]$ such that $f(a)A_0\subseteq B_1A_1$.
In particular, $f(a)\in B_1A_1$ and since $AA_1$ is nil, it is
done. For the second statement notice that $B_0$ is finitely
generated. Since $A$ is a finitely generated $A_0$-module, it
implies that $A$ is a finitely generated $B_0$-module. In
particular, $B_1$ is a finitely generated $B_0$-module.
\begin{pr}
Assume that $K[X]$ is a finitely generated $R$-module and $K[X]_0$
is a finitely generated algebra. Then $\tilde{X/G}\simeq SSp \ R$,
provided $G$ acts freely on $X$. Besides, $K[X]$ is a projective
$R$-module.
\end{pr}
Proof. One has to {\it superize} \cite{dg}, III, 4.6. More
precisely, we prove that $\psi$ is a superalgebra isomorphism and
$K[X]$ is a projective $R$-module. By Lemma 1.4 one can replace
$R$ and $K[X]$ by $R_{\cal P}$ and $K[X]_{\cal P}$, where $\cal P$
is a prime ideal of $R$. In other words, one can assume that $R$
is local and $K[X]$ is semi-local. Since $\psi$ is a
$K[X]$-supermodule morphism, Lemma 1.3 and Lemma 3.1 infer that
there are elements $f_1,\ldots , f_{m+n}\in K[X]$, where $|f_i|=0,
1\leq i\leq m=\dim K[G]_0, |f_i|=1, m+1\leq i\leq m+n, n=\dim
K[G]_1$, such that $\tau_X(f_1),\ldots , \tau_X(f_{m+n})$ form a
basis of the free $K[X]$-supermodule $K[X]\otimes K[G]$. Let $V$
be a superspace of superdimension $(m, n)$ with a basis
$v_1,\ldots , v_{m+n}$ such that $|v_i|=|f_i|, 1\leq i\leq m+n$.
Tensoring by $V$ the bottom line of the basic diagram we obtain a
diagram
$$
\begin{array}{ccccc}
K[X]\otimes K[G]^{\otimes 2} &
\begin{array}{c}{\stackrel{\delta'_0}{\leftarrow}} \\
{\stackrel{\delta'_1}{\leftarrow}}\end{array} & K[X]\otimes K[G]
&\stackrel{\tau_X}{\leftarrow} & K[X] \\
u_2 \uparrow & & u_1 \uparrow & & u_0 \uparrow \\
V\otimes K[X]\otimes K[G] & \begin{array}{c}{\stackrel{V\otimes\tau_X}{\leftarrow}} \\
{\stackrel{V\otimes i_X}{\leftarrow}}\end{array} & V\otimes K[X]
&\stackrel{i}{\leftarrow} & V\otimes R
\end{array},
$$
where $u_0(v_i\otimes r)=f_i r, u_1(v_i\otimes
f)=\tau_X(f_i)i_X(f), u_2(v_i\otimes
t)=\delta'_0(\tau_X(f_i))\delta'_2(t)$. By definition, $u_1$ is an
isomorphism of $K[X]$-supermodules. As in \cite{dg} we conclude
that $u_2$ is an isomorphism (of superspaces) and therefore, $u_0$
is. In particular, $K[X]$ is a free $R$-supermodule and the
elements $f_i$ form its basis. Returning to the general case, by
Lemma 1.5 from \cite{z} we obtain that $\psi$ is an isomorphism
and $K[X]$ is a projective $R$-module by \cite{f}, Theorem A.2.4.
By Lemma 1.5 (see also \cite{bur}, I, \S 2, Proposition 1) $K[X]$
is a faithfully flat (left and right) $R$-module. Proposition 4.2,
\cite{z}, concludes the proof.

Let a group $K$-sheaf $G$ acts on a $K$-sheaf $X$ freely. If $A,
B\in SAlg_K$ and $B$ is a fppf covering of $A$, then we denote
$B\succeq A$. Notice that $\succeq$ is a (partial) direct order.
If $X$ is a $K$-functor, the kernel
of maps $X(A)\begin{array}{c} \stackrel{X(i_1)}{\to} \\
\stackrel{X(i_2)}{\to}\end{array} X(B\otimes_A B)$, where
$i_1(a)=a\otimes 1, i_2(a)=1\otimes a, a\in A$, is denoted by
$X(B, A)$ (see \cite{z, jan} for more definitions and notations).
\begin{pr}
Let $N$ be a normal group $K$-subfunctor of $G$. Then the group
$K$-sheaf $\tilde{G/N}$ acts freely on $Y=\tilde{X/N}$ and
$\tilde{Y/H}\simeq\tilde{X/G}$.
\end{pr}
Proof. Denote the "naive" factors $$A\to G(A)/N(A), A\to
X(A)/N(A),
\\ A\in SAlg_K,$$ by $H_{(n)}$ and $Y_{(n)}$ correspondingly.
Consider $h\in H(A), y\in Y(A)$. There is a fppf-covering
$B\succeq A$ such that $h'=H(\iota^B_A)(g)\in H_{(n)}(B, A),
y'=Y(\iota^B_A)(y)\in Y_{(n)}(B, A)$. By the normality of $N$, the
group functor $H_{(n)}$ acts canonically on $Y_{(n)}$. In
particular, $y'h'\in Y_{(n)}(B, A)\subseteq Y(B, A)$. Since $Y$ is
a sheaf, one can define $yh=Y(\iota^B_A)^{-1}(y'h')\in Y(A)$. This
definition does not depend on the choice of $B$. In fact,  let $C$
be another fppf-covering of $A$. Then $D=B\otimes_A C\succeq B, C$
(see \cite{z, jan}). Thus $D\succeq A$. Set $$H(\iota^C_A)(h)=h'',
Y(\iota^C_A)(y)=y'', H(\iota^D_A)(h)=h''', Y(\iota^D_A)(y)=y'''
.$$ We have $$H(\iota^D_B)(h')=H(\iota^D_C)(h'')=h''',
Y(\iota^D_B)(y')=Y(\iota^D_C)(y'')=y''' .$$ It follows that
$Y(\iota^D_B)(y'h')=Y(\iota^D_C)(y''h'')=y'''h'''$. On the other
hand, all morphisms $Y(\iota^?_?)$ are mono and therefore,
$$Y(\iota^D_A)^{-1}(y'''h''')=Y(\iota^D_B)^{-1}(y'h')=Y(\iota^D_C)^{-1}(y''h'').$$
Similarly, one can prove that $H$ acts on $Y$ freely. To prove
that the above action is functorial on the argument $A\in SAlg_K$
one can mimic the proof of Lemma 2.3 from \cite{z}.

Finally, let $\rho : X\to Z$ be a $K$-sheaf morphism such that
$\rho(A)(xg)=\rho(A)(x)$ for all $x\in X(A), g\in G(A), A\in
SAlg_K$. There is a unique morphism $\alpha : Y\to Z$ satisfying
$\rho=\pi\alpha$, where $\pi : X\to Y$ is the canonical
factor-morphism. More precisely, for any $y\in Y(A)$ and for a
fppf covering $B\succeq A$ such that $Y(\iota^B_A)(y)=xN(B), x\in
X(B)$, we set $\alpha(A)(y)=Z(\iota^B_A)^{-1}(\rho(B)(x))$
\cite{jan, z}. Comparing with the definition of the $H$-action on
$Y$ we see that $\alpha$ is constant on $H$-orbits. In particular,
there is a unique morphism $\beta : \tilde{Y/H}\to Z$ such that
$\beta\pi'=\alpha$, where $\pi' : Y\to\tilde{Y/H}$ is the
corresponding factor-morphism. In other words, morphism $\pi'\pi :
X\to\tilde{Y/H}$ is the required factor-morphism. Theorem is
proved.
\begin{rem}
The same statement can be proved for dur $K$-sheafs.
\end{rem}
\begin{lm}
An (not necessary finite) algebraic supergroup $G$ acts freely on
an affine superscheme $X$ iff the ideal $J$ of $K[X]\otimes K[G]$,
generated by the elements $\tau_X(f)-f\otimes 1, f\in K[X]$,
contains $1\otimes {\cal M}$.
\end{lm}
Proof. If $1\otimes {\cal M}$ is not contained in $J$, then set
$A=K[X]\otimes K[G]/J$ and define  $$\alpha(f)=f\otimes 1 +J,
g(h)=1\otimes h +J, f\in K[X], h\in K[G].$$ It is obvious that
$g\in Stab_{G(A)}(\alpha)\setminus 1$. Conversely, if $g\in
Stab_{G(A)}(\alpha), \alpha\in X(A)$, then $\alpha\bar{\otimes}
g(\tau_X(f)-f\otimes 1)=0$ for any $f\in K[X]$, where
$\alpha\bar{\otimes} g(f\otimes h)=\alpha(f) g(h)$. The inclusion
$1\otimes {\cal M}\subseteq J$ implies $g({\cal M})=0$.

Now, everything is prepared to prove the main theorem. At first,
assume that $K[X]$ is finitely generated. Using induction on $|G|$
we prove that $K[X]$ is a finitely generated $K[X]^G$-module and
then apply Proposition 3.1 . If $G$ has a proper normal
supersubgroup $N$, then by the inductive hypothesis $K[X]$ is a
finitely generated $K[X]^N$-module. By Proposition 3.1
$\tilde{X/N}\simeq SSp \ K[X]^N$ and by Proposition 3.2
$\tilde{G/N}$ acts on $\tilde{X/N}$ freely. Since by Lemma 3.2
$K[X]^N$ is finitely generated, again the inductive hypothesis
infers that $K[X]^N$ is a finitely generated
$K[X]^G=(K[X]^N)^{\tilde{G/N}}$-module. Thus $K[X]$ is a finitely
generated $K[X]^G$-module. So, it remains to prove that $K[X]$ is
a finitely generated $K[X]^G$-module, whenever $G$ has not proper
normal supersubgroups. In particular, $G$ is either connected or
purely even and etale (cf. \cite{z}). Assume that $G$ is connected
and $G\neq 1$.
\begin{pr}
If $char K=0$, then $G\simeq G_a^-$.
\end{pr}
Proof. As it was noticed in \cite{z}, ${\cal M}=r+K[G]e$, where
$r$ is the radical of $K[G]$ and $e$ is the sum of primitive
idempotents belonging to ${\cal M}$. Besides,
$I_{G^{(0)}}=\bigcap_{t\geq 0}{\cal M}^t=K[G]e$. Since
$G=G^{(0)}$, it follows that $e=0$ and ${\cal M}=r=K[G]K[G]_1$. In
particular, ${\cal M}/{\cal M}^2$ is purely odd that implies
$Lie(G)_0=0$ and $Lie(G)_1^2=0$. In other words, $Lie(G)$ is
abelian and as in \cite{z} we conclude that $G$ is abelian. By
Lemma 9.5, \cite{z}, for any finite-dimensional $G$-supermodule
$V$ the equality $V^G= V^{Lie(G)}$ holds. In fact, $V^G$ is
naturally identified with $Hom_G(K, V)$, where $K$ is regarded as
one-dimensional trivial $G$-supermodule. Identify $Lie(G)$ with an
odd abelian supersubalgebra $L$ of $gl(V)$. Then for all $x, y\in
L$ we have $xy=-yx$. Let $A$ be an associative subalgebra of
$End_K(V)$ without unit, generated by $L$. It is clear that
$A^{\dim L +1}=0$ and by Engel's theorem there is a vector $v\in
V$ such that $Av=0$. In particular, $G$ is unipotent. Proposition
2.2 concludes the proof.
\begin{rem}
Proposition 3.3 infers that over an algebraically closed field of
characteristic zero, any finite supergroup is an extension of
abelian unipotent supersubgroup by an even etale group. It seems
to be very likely that such extension have to be split (for the
classical case see Theorem 3.3 from \cite{ko}).  We hope to check
all details in a next article and to get rid of the assumption
about the ground field to be algebraically closed.
\end{rem}

Let $G=G^-_a$. Remind that $$K[G]=K[t], |t|=1,
\delta_G(t)=t\otimes 1 +1\otimes t, \epsilon_G(t)=0, s_G(t)=-t .$$
\begin{lm}
A $G$-supermodule structure on a superspace $V$  is uniquely
defined by an odd (locally finite) endomorphism $\phi : V\to V,
\phi^2=0$. Precisely, $\tau_V(v)=v\otimes 1 +\phi(v)\otimes t$ and
therefore, $V^G=\ker\phi$.
\end{lm}
Proof. Straightforward calculations.

By Lemma 3.4 $\tau_X(f)=f\otimes 1 +\phi(f)\otimes t, f\in K[X]$,
where $\phi\in End_K(K[X])_1$ and $\phi^2=0$. Since $\tau_X$ is a
superalgebra morphism, we obtain that $\phi(f_1 f_2)=f_1\phi(f_2)
+(-1)^{|f_2|}\phi(f_1)f_2, f_1, f_2\in K[X]$. In other words,
$\phi$ is a right (odd) superderivation.
\begin{lm}
The supergroup $G$ acts on $X$ freely iff there is $f\in K[X]_1$
such that $\phi(f)\in K[X]^*$.
\end{lm}
Proof. By Lemma 3.3 $G$ acts on $X$ freely iff there are
$h_1,\ldots , h_n\in K[X]_1$ and $f_1,\ldots , f_n\in K[X]_0$ such
that $\sum_{1\leq i\leq n}f_i\phi(h_i)=1$. Thus
$$\phi(\sum_{1\leq i\leq n} f_i h_i)= 1- \sum_{1\leq i\leq
n}\phi(f_1)h_i\in K[X]^* .$$
\begin{lm}
The superalgebra $K[X]$ is a free $R$-supermodule of rank $2$.
\end{lm}
Proof. Consider $g\in K[X]_1$ such that $g=\phi(f)\in K[X]^*$. Set
$z=fg^{-1}$. Since $\phi(z)=1$, for any $h\in K[X]$ we have
$\phi(hz)=h-\phi(h)z$. Thus $K[X]=R + Rz$. If $h\in R\bigcap Rz$,
then $h=rz, r\in R,$ and therefore, $0=\phi(h)=r$.
\begin{exam}(see \cite{z}, section 10)
Consider a $G=G^-_a$-supermodule $V$ with a basis $v_1, v_2,
|v_1|=0, |v_2|=1,$ such that $\phi(v_1)=v_2, \phi(v_2)=0$ in the
above notations. The symmetric superalgebra $S(V)$ has the induced
$G$-supermodule structure by
$$\tau_{S(V)}(v_1^r)=v_1^r\otimes 1 +rv_1^{r-1}v_2\otimes t, \
\tau_{S(V)}(v_1^{r-1}v_2)=v_1^{r-1}v_2\otimes 1, r\geq 0.$$ By
Lemma 3.3 the induced $G$-action on $X=SSp \ S(V)$ is not free.
Moreover, $K[X]=S(V)$ is nor finitely generated $R$-module neither
integral over $B$, provided $char K=0$. In fact, the superalgebra
$R=K\bigoplus (\bigoplus_{r\geq 1}K v_1^{r-1}v_2)$ is not finitely
generated. In \cite{z} it was also proved that $K[X]$ is not any
flat $R$-module. Notice that the final conclusion in \cite{z} is
not completely correct. Indeed, Proposition 4.1 holds for free
actions which is not the case.
\end{exam}
Let $char K=p> 0$. An algebraic  supergroup $H$ is called {\it
infinitesimal supergroup of hight} 1 iff $h^p=0$ for any $h\in
{\cal M}$.
\begin{lm}
If $char K=p> 0$ and $G$ is connected, then $G$ is infinitesimal
supergroup of hight 1. In particular, $K[X]$ is a finitely
generated $K[X]^G$-module.
\end{lm}
Proof. As above, ${\cal M}=\ker\epsilon_G=r$. We have a series
$1\leq G_1\leq G_2\leq\ldots\leq G$, where each $G_n$ is a $n$-th
infinitesimal supersubgroup (cf. \cite{z, jan}). Since any $G_n$
is a normal supersubgroup of $G$, we have either $G_1=1$ and
$K[G]=F(K[G])=\{f^p |f\in K[G]\}$ (see the notice before Lemma
8.2, \cite{z}), or $G_1=G$. The equality $K[G]=F(K[G])$ implies
$r=F(r)$. The nilpotence of $r$ infers $r=0$ and $G=1$. The last
case $G=G_1$ is equivalent to $f^p=0$ for any $f\in {\cal M}$.
Finally, for any $f\in K[X]$ we have $\tau_X(f)=f\otimes 1 +\sum
f_1\otimes h_2$, where each $h_2$ belongs to $\cal M$. Thus
$\tau_X(f^p)=f^p\otimes 1$, that is $f^p\in K[X]^G$. Lemma 3.2
concludes the proof.

Now, let $G$ be even and etale. It is well known that $K[G]\simeq
(K\Gamma)^*$, where $\Gamma$ is a finite group and $K\Gamma$ is
its group algebra, endowed with Hopf algebra structure by
$\delta_{K\Gamma}(\gamma)=\gamma\otimes\gamma,
s_{K\Gamma}(\gamma)=\gamma^{-1}, \gamma\in\Gamma$ (see \cite{jan},
part I (8.5, 8.21) and \cite{w}, 2.3, 6.4, or see \cite{dg}, II,
\S 5, 2.4). Therefore, $K[G]$ is generated by the idempotents
$e_{\gamma}$, such that
$e_{\gamma}(\gamma')=\delta_{\gamma,\gamma'}$ and
$$\epsilon(e_{\gamma})=\delta_{\gamma, 1}, \delta_G(e_{\gamma})=\sum_{\gamma'\in\Gamma}e_{\gamma'}\otimes e_{\gamma'^{-1}\gamma},
s_G(e_{\gamma})=e_{\gamma^{-1}}, \gamma,\gamma'\in\Gamma.$$ A
vector superspace $V$ is called $\Gamma$-{\it supermodule} iff it
is a $\Gamma$-module and any $\gamma\in\Gamma$ acts on $V$ as an
even operator. The category of $\Gamma$-supermodules with even
morphisms is denoted by $\Gamma-smod$. If $V\in\Gamma-smod$, then
it has a $G$-supermodule structure by
$$\tau_V(v)=\sum_{\gamma\in\Gamma}\gamma v\otimes e_{\gamma}, v\in
V.$$ This correspondence defines an equivalence of categories. In
particular, $G$ acts on an affine superscheme $X$ iff $K[X]$ is a
$\Gamma$-supermodule and any $\gamma\in\Gamma$ acts as a
superalgebra automorphism. Since $$K[X]^G=K[X]^{\Gamma}=
K[X]_0^{\Gamma}\bigoplus K[X]_1^{\Gamma}$$ this case is also done.

It remains to consider the case when $K[X]$ is not finitely
generated. Since any $K[G]$-supercomodule is locally finite, the
superalgebra $K[X]$ is a direct union of its finitely generated
subalgebras $B_i, i\in I,$ such that each $B_i$ is a
$G$-supersubmodule of $K[X]$. In other words, $G$ acts on any $SSp
\ B_i$ and the canonical morphism $SSp \ B_i\to X$ commutes with
this action. Since $\cal M$ is finite-dimensional, by Lemma 3.3
one can assume that $G$ acts freely on each $SSp B_i$. By the
above, for any $i\in I$ the superalgebra $B_i$ is a faithfully
flat (left and right) $R_i=B_i^G$-module and the canonical
morphism $B_i\otimes_{R_i} B_i\to B_i\otimes K[G]$ is an
isomorphism. Thus $K[X]=\lim\limits_{\rightarrow} B_i$ is a
faithfully flat (left and right) $R=\lim\limits_{\rightarrow}
R_i$-module (cf. Lemma 7.1, III, \S 3, \cite{dg}) and
$$K[X]\otimes_R K[X]=\lim\limits_{\rightarrow} B_i\otimes_{R_i}
B_i\simeq \lim\limits_{\rightarrow} B_i\otimes K[G]=K[X]\otimes
K[G].$$ Use Proposition 11, \cite{bur}, I, \S 3, and the above
isomorphism (of $K[X]$-modules) one can conclude that $K[X]$ is
finitely presented. Exercise 15 from \cite{bur}, I, \S 2, implies
that $K[X]$ is also a projective $R$-module. Remark 1.1 and
Proposition 4.2, \cite{z}, infer  $\tilde{X/G}\simeq SSp \ R$.

\begin{center}
\bf Acknowledgements
\end{center}
This work was partially supported by RFFI 07-01-00392.

\end{document}